\documentclass[11pt]{article} % use larger type; default would be 10pt

\usepackage[utf8]{inputenc} % set input encoding (not needed with XeLaTeX)

\usepackage{geometry} % to change the page dimensions
\usepackage{graphicx} % support the \includegraphics command and options

\usepackage{booktabs} % for much better looking tables
\usepackage{array} % for better arrays (eg matrices) in maths
\usepackage{paralist} % very flexible & customisable lists (eg. enumerate/itemize, etc.)
\usepackage{verbatim} % adds environment for commenting out blocks of text & for better verbatim
\usepackage{subfig} % make it possible to include more than one captioned figure/table in a single float
\usepackage{amsmath}
\usepackage{amssymb}
\newtheorem{theorem}{Theorem}

\usepackage{fancyhdr} % This should be set AFTER setting up the page geometry
\usepackage{sectsty}
\allsectionsfont{\sffamily\mdseries\upshape} % (See the fntguide.pdf for font help)
\usepackage[nottoc,notlof,notlot]{tocbibind} % Put the bibliography in the ToC
\usepackage[titles,subfigure]{tocloft} % Alter the style of the Table of Contents

\title{Admissible sets of integers in imaginary quadratic fields of class number 1}
\author{Daniel Shiu}
\begin{document}
\maketitle

\section{Introduction}
Let $\mathcal O_K$ be the ring of integers of a number field $K$. A finite subset $\mathcal H=\{\eta_1,\ldots,\eta_k\}\subset\mathcal O_K$ is said to be \emph{admissible} if 
$$\forall \mathfrak p\subseteq\mathcal O_K\quad \mathcal H/\mathfrak p\neq\mathcal O_K/\mathfrak p$$
i.e. for every prime ideal $\mathfrak p$ in $\mathcal O_K$, the reductions of the $\eta_i$ modulo $\mathfrak p$ do not represent every residue class. We note that for prime ideals $\mathfrak p:\mathcal N\mathfrak p>k$ this condition is automatically satisfied as the number of residue classes is greater than the number of elements in the set. The diameter of a finite subset $S$ of $\mathcal O_K$  is defined to be
\[\mathrm{diam}(S):= \max_{\eta_i,\eta_j\in S:\eta_i\neq\eta_j}|\mathcal N(\eta_i-\eta_j)|\]
where $\mathcal N$ denotes the norm down to the rational field. This paper computationally investigates admissible subsets of rings of integers in imaginary quadratic fields of class number 1 (such as the Gaussian integers or the Eisenstein integers) that are large relative to their diameter. The motivation is twofold. 

Firstly, there is the improvement of results analogous to the question of an upper bound for the least gap between large rational primes, which has been one of the great number theoretic successes of the twenty-first century. Starting from the work of Goldston-Pintz-Yildirim \cite{goldston2009primes}, Zhang \cite{zhang2014bounded}, Maynard \cite{maynard2015small} and Polymath \cite{polymath2014variants} have been able to show that for any $m$, there is a bound $k(m)$ such that any admissible set of at least $k(m)$ elements has infinitely many translates for each of which at least $m$ translated elements are primes. The case that attracts most attention is the case $m=2$ for which the Polymath project demonstrated the bound $k(2)\le 50$. In general if we can construct admissible sets of rational integers of size at least $k(2)$ and diameter at most $D$ we may conclude
\[\liminf_{n\to\infty}p_{n+1}-p_n\le D\]
In the rational case the best admissible sets of rational integers of size 50 have diameter 246, which was the construction used by the Polymath project. Subsequent work by Stadlman \cite{stadlmann2026bounded} has claimed the improvement $k(2)\le 49$ for which the best admissible set has diameter 240. Further work by Song et al \cite{song2026bounded} has claimed the improvement $k(2)\le 48$ for which there the best admissible set has diameter 236.

The result has been generalised to number fields by Castillo et al \cite{castillo2015bounded}, and Vatwani \cite{vatwani2017bounded} considered the particular case of the Gaussian integers (using methods that extend to other quadratic fields of class number 1) for which admissible sets of size $k(2)=50$ can be shown to have translates that capture at least two primes. It seems likely that the methods of Stadlman and Song et al can be similarly applied so that $k(2)\le 48$ for such fields. Vatwani did not optimise the choice of admissible sets, instead electing to use the same admissible set of rational integers as the Polymath project. Using the norm for the measure of our diameter this shows that
\begin{equation}\liminf_{\rho,\rho'{\ \rm prime:\ }\mathcal N(\rho)\to\infty}\mathcal N(\rho-\rho')\le 246^2\label{eq:vatwani}\end{equation}
with the same upper bound when we restrict to $\rho-\rho'\in\mathbb Z $. The constructions in this paper improve these bounds by improving the admissible sets, both in the case where we consider admissible sets like those of Vatwani, all of whose elements differ by rational integers, or in the general case.

In particular then, we offer the following improvement/extension of Vatwani's theorem 1.4/theorem 13.1 \cite{vatwani2017bounded}
\begin{theorem}\label{th:line}
There are infinitely many pairs of rational primes of the form $p_1=a^2 +b^2$ and $p_2 =(a+h)^2 + b^2$, with $a,b,h\in\mathbb Z$, such that $0 <|h| \le 136$.

Similarly, there are infinitely many pairs of rational primes of the form $p_1=a^2 +2b^2$ and $p_2 =(a+h)^2 +2b^2$, with $a,b,h\in\mathbb Z$, such that $0 <|h| \le 168$.

Further, there are infinitely many pairs of rational primes of the form $p_1=a^2 - ab +db^2$ and $p_2 =(a+h)^2 - (a+h)b+db^2$, with $a,b,h\in\mathbb Z$, $d=1, 2, 3, 5, 11, 17, 41$ and $0<|h|\le  93, 126, 100, 80, 59, 55, 50$ respectively.
\end{theorem}
These results can all be simply improved if the Stadlman and Song work transfers, and the corresponding improvements are given at the end of section \ref{sec:lines}.

Limiting our admissible sets to those whose elements differ by a rational integer is a significant constraint, and the two-dimensional structure of our rings of integers suggest instead considering sets of integers contained in a region of bounded diameter with large area, such as a disc. It turns out that methods for searching for rational admissible sets do not automatically transfer to this setting, but improvement over the linear structure is relatively straightforward to achieve. In particular then, we improve/extend Vatwani's theorem 12.3 (\ref{eq:vatwani}) as follows:
\begin{theorem}\label{th:disc}
There are infinitely many pairs of distinct Gaussian primes $\pi_1$, $\pi_2$ such that ${\mathcal N}(\pi_1-\pi_2) \le 218$. 

Furthermore for each of the rings of integers for the fields $\mathbb Q(\sqrt{-\Delta})$ for 
$$\Delta=3, 7, 8, 11, 19, 43, 67, 163$$ there are infinitely many pairs of prime elements $\pi_1$ and $\pi_2$ such that 
$${\mathcal N}(\pi_1-\pi_2)<147, 452, 454, 441, 397, 341, 400, 503$$ 
respectively.
\end{theorem}

A second motivation is that there is the opportunity to disprove a variation of the second Hardy-Littlewood conjecture \cite{hardy1923some}, which posits for the rational prime counting function $\pi(x)$:
\begin{equation}\text{?`}\quad\forall x,y\in\mathbb Z, x,y \ge 2 \Rightarrow\pi(x+y)\le\pi(x)+\pi(y)\quad\text{?`}.\label{eq:hl2}\end{equation}
Ironically, Hensley and Richards \cite{hensley1973incompatibility} were able to show that this is incompatible with the first Hardy-Littlewood conjecture \cite{hardy1923some} which posits (inter alia) that there are infinitely many shifts of an admissible set such that all elements are prime. Hensley and Richards did this by exhibiting an admissible set contained in an interval $[1,y_0]$ which contains more that $\pi(y_0)$ elements. It follows then that for any shift $x$ for which the first Hardy-Littlewood conjecture captures the full set of primes we have
\[\pi(x+y_0)-\pi(x)>\pi(y_0).\]
One generalisation of the second Hardy-Littlewood conjecture to imaginary quadratic number fields of class number 1 would be to posit that that the primes in the field are more concentrated around the origin. If we write $\pi(\xi,A)$ for the count of prime elements (counting associates with multiplicity) in a disc $\{\eta:|\mathcal N(\xi-\eta)|\le A\}$ centred on the algebraic integer $\xi$, a plausible analgoy would be the conjecture that
\begin{equation}\text{?`}\quad\pi(0,A)\ge\pi(\xi,A)\quad\text{?`}\label{eq:ghl2}\end{equation}
for all $\xi\in\mathcal O_K$ and $A$ such that $\pi(0,A)>0$. In the spirit of Hensley and Richards we might hope to exhibit an admissible set with cardinality greater than $\pi(0,A_0)$ contained in a disc $\mathcal N(\eta)\le A_0$, for some $A_0$. This would demonstrate an incompatibility with an analogy to the first Hardy-Littlewood conjecture. Going further, depending on the size of $A_0$, we might even hope to find a relatively small shift of the set which captures all of its primes, giving an explicit counter-example. The ``lucky number'' effect noted by Euler shows that our choices of fields are sparse in small primes as the discriminant grows, which gives us hope that $\pi(0,A)$ is sufficiently smaller or that our admissible sets are denser than in the rational case and that we might find explicit counter-examples with modest computational effort.

Note that (\ref{eq:ghl2}) is a slightly weaker claim than (\ref{eq:hl2}) as we count associates of primes with multiplicity. We adjust (\ref{eq:hl2}) to a comparable, strictly weaker conjecture due to Erd\"os \cite{erdos1981problems} about the primes near 0 where we consider primes (counting $p$ and $-p$ separately) in the interval $[-y/2,y/2]$ with $y\ge 4$. The weaker conjecture then becomes
\begin{equation}\text{?`}\quad\forall x,y\in\mathbb Z, x,y \ge 4 \Rightarrow\pi(x+y)\le\pi(x)+2\pi(y/2)\quad\text{?`}\label{eq:chl2}\end{equation}
A counterexample now requires the stronger property
\[\pi(x+y_0)-\pi(x)>2\pi(y_0/2),\]
but suitable admissible sets of rational integers exist that again show that this is incompatible with the first Hardy-Littlewood conjecture per the work of Clark and Jarvis \cite{clark2001dense}.

The methods for constructing admissible sets in this paper quickly lead to explicit counter-examples in most cases, though a frustrating omission is the case $\mathbb Q(\sqrt{-11})$ where we have failed to construct a sufficiently dense admissible set, let alone an instance where a shift fully captures all of its primes. Nevertheless we present:

\begin{theorem}
Let $\mathcal O_K$ denote the rings of integers of a number field and $\pi(\xi,A)$ for the count of prime elements in the disc $\{\eta:|\mathcal N(\eta-\xi)|\le A\}$ centred on $\xi\in\mathcal O_K$. Then the conjecture
\[\pi(0,A)>0\Rightarrow \forall\xi\in\mathcal O_K\quad  \pi(0,A)\ge\pi(\xi,A)\]
does \emph{not} hold in the fields $K=\mathbb Q(\sqrt{-\Delta})$ for
$$\Delta=3,4,7, 8, 19, 43, 67,163.$$
\end{theorem}

The question of finding dense admissible sets is related to, but not identical to the ``Gaussian zoo'' question \cite{renze2001gaussian} where the search is for large admissible connected sets where pairs of Gaussian primes are connected if their difference has norm 2.

ChatGPT was used to source some of the background references to this paper and to offer corrections to its first draft, but otherwise no Large Language Model was used in this work.

\section{Short line segments containing admissible sets\label{sec:lines}}
The admissible set by Vatwani \cite{vatwani2017bounded} is the same used in \cite{polymath2014variants} which was found exhaustively after larger admissible sets were constructed using variations on the ``Schinzel sieve'' \cite{schinzel1961}. The general non-exhaustive method is to consider a small interval of candidates containing 0 (sets remain admissible under arbitrary translation, so that this choice is without loss of generality) and to specify a residue class mod $p$ for each prime less than the target size of the admissible set. We delete integers from our list of candidates if they lie in any of the specified residue classes; we can then choose an admissible set of the target lower bound from the remaining candidates provided that there are at least as many survivors as our target bound. The goal then is to choose the residue classes to maximise the number of survivors. Schinzel \cite{schinzel1961} observed that an asymptotically effective way of choosing residue classes is to choose the residue class 1 for primes up to some bound $z\approx \sqrt{\log\log k}$ and then specify the residue class 0 for large primes. For $k<5.14\times 10^{23}$ we have $\sqrt{\log\log k}<2$, but nevertheless, from this starting point, Polymath found good admissible sets with a somewhat \emph{ad hoc} method. They choose the residue class 1 mod 2, then choose the residue classes 0 mod $p$ for odd primes less than roughly $\sqrt{k\log k}$, then greedily choose the residue class mod $p$ that eliminates fewest survivors for large $p$. After sieving out the set of residue classes, further small improvements were achieved by small perturbations of the interval endpoints and ``swapping'' small sets of residue classes of the larger primes. Eventually, the target size came within the bounds of the exhaustive search by Engelsma, who lists an optimal admissible set of size 50:
\begin{verbatim}
  0   4   6  16  30  34  36  46  48  58  60  64  70  78  84  88  90  94 
100 106 108 114 118 126 130 136 144 148 150 156 160 168 174 178 184 190
196 198 204 210 214 216 220 226 228 234 238 240 244 246
\end{verbatim}

Admissible sets of rational integers are also admissible sets of algebraic integers in any number field. To see this note that we simply have to specify a residue class modulo each prime ideal of norm less than the size of our set, and that specifying the residue class corresponding to the residue class for the rational prime lying above the ideal in our rational set will preserve all of our survivors. Hence Vatwani's choice of admissible set is valid. It is however suboptimal. 

Consider the Gaussian integers. If we elect to build an admissible set from a collection of rational integers, we are allowed, say, to specify the residue class $i\pmod 3$ this requires us to eliminate all candidates of the form $\{3a+(3b+1)i:a,b\in\mathbb Z\}$, none of which are rational integers. Similarly, for any rational prime that is $3\pmod 4$ we find that the inertia of the prime means that we can specify the residue class $i\pmod p$ and no rational integer is eliminated. For rational primes that are $1\pmod 4$, these split into pairs of ideals and no matter what residue class we pick every $p$th rational integer can represent that class. However, we can specify that the residue class for both ideals contains the same rational integers so that we lose a proportion roughly $1/p$ of our candidates rather than $2/p$. We deduce that we can sieve an interval of rational integers selecting a residue class mod 2 and a residue class for each rational prime $p:p<k,p\equiv 1\pmod 4$ and if we have more than $k$ survivors, we have an admissible set of Gaussian integers. 

For $k=50$, we need to sieve by the primes $\{2,5,13,17,29,37,41\}$ and an exhaustive search shows that the smallest diameter such set has rational diameter $D=136$ (hence Gaussian diameter $136^2$). In fact there are 272 subintervals of this length mod $2\times 5\times 13\times 17\times 29\times 37\times 41$, the first being $[68917,69053]$. For $k=49$ the least rational diameter is $D=133$ and for $k=48$, the least rational diameter is $D=131$. It is worth noting that following the Polymath method without small perturbations recovers these admissible sets, showing that the method is very effective at these sizes of $k$.

The same ``avoid inert primes; double-up split primes'' approach works for other imaginary quadratic fields of class number 1. Thus we have  admissible sets with 50 elements (resp. 49, 48), whose members differ by rational integers and which have diameter at most $D$ for the values of $D$ given in the table in Figure \ref{tab:line} for each such number field. We list the number fields in order of increasing absolute value of discriminant.

\begin{center}
\begin{figure}[htbp]
\centering
\begin{tabular}[h]{|c|c|c|c|c|}
\hline
Discriminant & Split primes & $D: k=50$ & $D: k=49$ & $D: k=48$\\
\hline
-3 & \{3,7,13,19,31,37,43\} & 93 & 91 & 88\\
-4 & \{2,5,13,17,29,37,41\} & 136 & 133 & 131\\
-7 & \{2,7,11,23,29,43\} & 126 & 122 & 120\\
-8 & \{2,3,11,17,19,41,43\} & 168 & 164 & 162\\
-11 & \{3,5,11,23,31,37,47\} & 100 & 98 & 96\\
-19 & \{5,7,11,17,19,23,43,47\} & 80 & 78 & 77\\
-43 & \{11,13,17,23,31,41,47\} & 59 & 58 & 57\\
-67 & \{17,19,23,29,37,47\} & 55 & 54 & 52\\
-163 & \{41,43,47\} & 50 & 49 & 48\\
\hline
\end{tabular}
\caption{Admissible sets in an interval results}
\label{tab:line}
\end{figure}
\end{center}

In particular, the $k=50$ column of the table in Figure \ref{tab:line} delivers Theorem \ref{th:line}, with other columns giving results in line with anticipated improvements to the $k(2)$ value. Specific examples of such admissible sets for $k=50$ are given in the appendix.

We note that multiplying an admissible set by an element of our ring of integers that is coprime to our sieving set also produces an admissible set and if the element is not coprime, matters may improve as we can remove prime divisors of the elements from our sieving set in the same way that we avoid inert primes. This computation therefore also delivers bounded norms of differences between primes in our fields where the difference is a rational integer multiple of a fixed element of the ring (scaling the bound by the norm of the fixed element in question).

\section{Discs containing admissible sets\label{sec:discs}}
A more interesting question arises when we consider general differences, so that we are considering the two-dimensional distribution of the integers in our ring. If we seek an admissible set of small diameter, a natural set of initial candidates will be the intersection of our ring of integers with a disc $\{\eta:\mathcal N(\eta)<B\}$, which we sieve with prime ideals $\mathfrak p:\mathcal N\mathfrak p<=k$ by specifying a residue class modulo each prime ideal. Survivors are those candidates which lie in none of the specified residue classes, and if we have more than $k(2)$ survivors, we can guarantee infinitely many pairs of primes the norm of whose difference is less than $4B$.

In order to anchor our expectations for the results of our investigations, we begin with a Zhang-like approach, where we specify the residue class 0 for all of the Gaussian prime ideals of norm less than 50. We list the Gaussian integers in order of increasing norm (tie-breaking on argument) and proceed down the list until we find the 50th survivor. This occurs with our 323rd candidate, which has norm 101 and lies in a set of (Gaussian) diameter 404.

Taking our lead from Schinzel, we could instead specify the residue class 1 for the prime ideal $\langle 1+i\rangle$ and the residue class 0 for the other prime ideals of norm less than 50. In this case our 50th survivor occurs when we sieve a disc with 489 candidates. Our largest survivor has norm 153 and lies is a set of (Gaussian) diameter 592. This gives a first indication that our methods may not transfer very well.

Following Polymath, we can specify residue class 1 for the prime ideal $\langle 1+i\rangle$, the residue class 0 for the other prime ideals of norm less than $\sqrt{50\log 50}\approx 13.99$, and the residue class that eliminates fewest survivors for the remaining prime ideals. In this case our 50th survivor occurs with our 250th candidate which has norm 81 and lies is a set of (Gaussian) diameter 298. The ``greedy'' sieving phase, eliminates only two survivors, which encourages us to make greater use of greedy strategies. Lowering the bound at which we being the greedy approach, seems to consistently improve the admissible sets found.

Without committing ourselves to the residue class zero for any prime (so that we specify residue class 1 for the prime ideal $\langle 1+i\rangle$ and select all other residue classes greedily) yields further improvement. We refer to this strategy as ``odd greedy''; it leads to our 50th survivor occurring with our 190th candidate which has norm 61 and lies is a set of (Gaussian) diameter 244.

Allowing even the residue class for $\langle 1+i\rangle$ to be chosen greedily, we have a strategy which we call ``pure greed'' that  leads to our 50th survivor occurring with our 182nd candidate which has norm 58 and lies is a set of (Gaussian) diameter 218.

\begin{center}
\begin{figure}[htbp]
\centering
\begin{tabular}{m{0.5\linewidth} m{0.5\linewidth}}
\includegraphics[height=5cm]{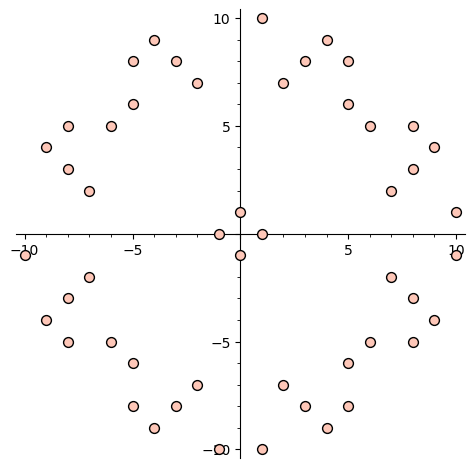} &
\includegraphics[height=5.5cm]{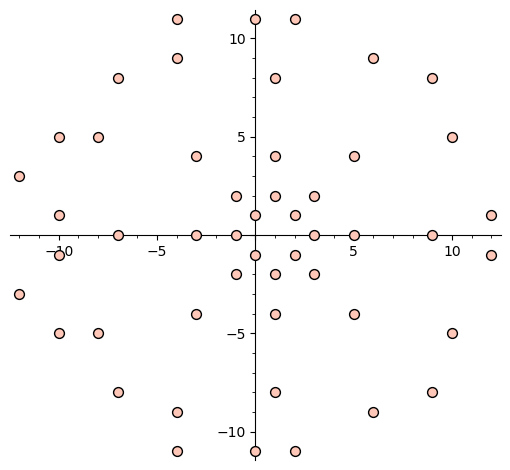} 
\end{tabular}
\caption{Zhang and Schinzel approaches}
\label{fig:sieve}
\end{figure}
\end{center}

\begin{center}
\begin{figure}[htbp]
\centering
\begin{tabular}{m{0.3\linewidth} m{0.3\linewidth} m{0.3\linewidth}}
\includegraphics[height=4.5cm]{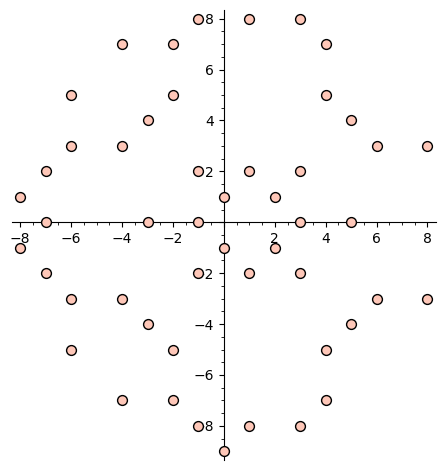} &
\includegraphics[height=3.5cm]{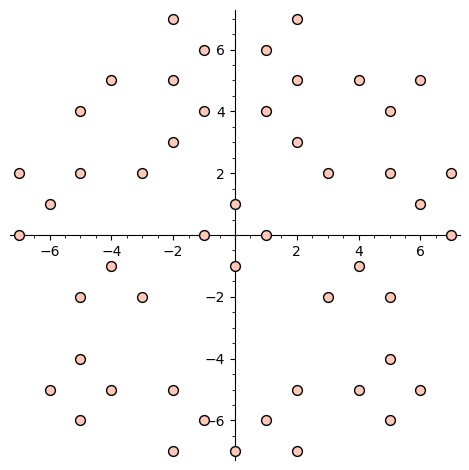} &
\includegraphics[height=3.5cm]{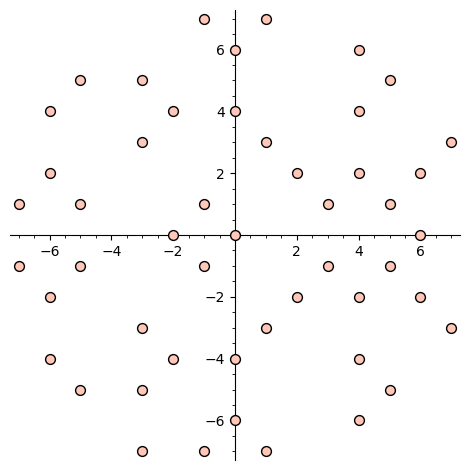} 
\end{tabular}
\caption{Polymath, odd greedy, and pure greedy approaches}
\label{fig:greed}
\end{figure}
\end{center}

It is instructive to visualise the layouts of the admissible sets produced by these five methods and they are pictured in five Argand diagrams in Figures \ref{fig:sieve} and \ref{fig:greed}. We see that the Zhang method (Figure \ref{fig:sieve} left) has four survivors in the unit locations, but then has a wide ``moat'' before permitting a ring of survivors corresponding to primes whose norms are greater than $k$. Schinzel's sieve (Figure \ref{fig:sieve} right) is intended to provide a handful more survivors in the centre while simultaneously adding density to the outer ring at the price of increasing the width of the moat. This can be viewed as attempting to structurally optimise some of the inclusion/exclusion of candidates. The Polymath sieve (Figure \ref{fig:greed} right) intends to sacrifice some of the survivors in the ring in the hope of generating more survivors in the moat area. We might interpret this as using Schinzel to structurally optimise the inclusions/exclusions involving only small primes, then using greed to unstructurally optimise the inclusions/exclusions that involve large primes. The odd greedy (Figure \ref{fig:greed} centre) and pure greedy (Figure \ref{fig:greed} right) approaches indicate that even when we do not build structure into our choices, the process still captures some helpful structure. On the other hand, it also suggests that the optimisation of the inclusion/exclusion effects is somehow different from the  analogue in the rational integers.

One plausible explanation for the difference between the rational and quadratic field cases is the greater non-uniformity in the representation of residue classes within our chosen region. In the rational case, our region is an interval of the real line of length, say, $L$, and the number of integer candidates eliminated by any single residue class modulo a prime $p$ is either $[L/p]$ or $[L/p]+1$. In the Gaussian case with a disc of area, say $A$, if we consider the residue class zero modulo some prime $\rho$, the candidates eliminated correspond to the integer points inside a disc of area $A/\mathcal N\rho$, which is Gauss's famous circle problem. Consequently, per the work of Hardy \cite{hardy1915expression} and Landau, the number of eliminations could be 
\begin{equation}\frac A{\mathcal N\rho}+\Omega_{\pm}((A\log(A/\mathcal N\rho)/\mathcal N\rho)^{1/4})\label{eq:circle_problem}\end{equation} 
with similar discrepancies for other residue classes. This suggests that even absent inclusion/exclusion discrepancy accumulations, there is potentially very large value for choosing which residue classes to eliminate based on their interaction with the boundary of our region. This will be particularly true of small primes where the $\Omega_{\pm}$ term in equation (\ref{eq:circle_problem}) is largest for a given fixed disc.

One might hope to simultaneously optimise the boundary variation and the inclusion/exclusion by performing a greedy approach selecting residues of composite moduli. We might take a product, $\mathfrak P$ of small prime ideals whose norm slightly greater than the number of candidates, then for each residue class modulo $\mathfrak P$, count the number of candidates whose difference from the residue class is coprime to $\mathfrak P$. We could then specify the residue class mod $\mathfrak P$ with the smallest count (which then determines the residue class specified for each prime ideal factor of $\mathfrak P$). This ``sieve by composite'' approach is more computationally burdensome than sieving by prime, though is feasible for moderately sized sets of candidates.

For the Gaussian integers with a target size of $k=50$, the ``greedy composites'' approach matched the greedy prime approach. However, for other fields an improvement was obtained.

We further note that the two-dimensional nature of our investigation offers another degree of freedom not granted to the rational integer case. We chose our initial set of candidates as those integers contained in a disc centred on the origin, which is natural as it has a perimeter of constant diameter providing an upper bound for the diameter of our final admissible set. In the rational integers, there is really only one such choice of an interval (up to fencepost effects), but in a two-dimensional setting there is more variation possible with our choice of the centre of the disc (or more ambitiously we might choose a different shape of constant diameter such as a  Reuleaux triangle). Plotting the Minkowski embedding of the integers in our fields, we see translates of the rational integer line that grow in separation as the discriminant of our field grows. Choosing a disc centred on the origin concentrates our candidate on a single line as the disc radius shrinks relative to the discriminant, however as a first step, it is worth investigating whether centring between two lines significantly changes our outcomes. For each of the fields we therefore also construct starting sets centred on $(1+\alpha)/2$ where $\alpha$ is the other canonical integral basis element. We can experiment with both the ``greedy prime'' and ``greedy composite'' methods (the Zhang, Schinzel, and Polymath methods become ambiguous). In the case of the Gaussians (centre $(1+i)/2$) with target size $k=50$, this choice of initial set was worse than the origin centred choice, but for other fields there was some improvement. We summarise the size of the smallest initial set and diameter of the admissible set found for each of the four possible approaches for each field in the table in Figure \ref{tab:disc}
\begin{center}
\begin{figure}[htbp]
\centering
\begin{tabular}[h]{|c||c|c|c|c|}
\hline
Discriminant & Prime, 0 & Comp., 0 & Prime, $(1+\alpha)/2$ & Comp. $(1+\alpha)/2$\\
\hline
-3 & (145, 151) & (147,156) & (143,151) & (140,{\bf147)} \\ 
-4 & (182,{\bf 218)} & (182,{\bf 218)} & (190,234) & (195,234) \\
-7 & (281,464) & (282,{\bf 452)} & (302,484) & (287,464) \\
-8 & (262,450)  & (262,450) & (258,{\bf 454)} & (258,{\bf 454)}\\
-11 & (229,484) & (219,452) & (221,444) & (217,{\bf 441)} \\
-19 & (160,448) & (159,425) & (156,{\bf 397)} & (152,400) \\
-43 & (85,368) & (87,368) & (89,{\bf 341}) & (89,{\bf 341})\\
-67 & (81,{\bf 400}) & (81,{\bf 400}) & (80,421) & (80,425) \\
-163 & (66,576) & (66,576) & (71,523) & (70,{\bf 503})\\
\hline
\end{tabular}
\caption{Admissible sets inside discs results}
\label{tab:disc}
\end{figure}
\end{center}

\begin{center}
\begin{figure}[htbp]
\centering
\begin{tabular}{m{0.33\linewidth} m{0.33\linewidth} m{0.33\linewidth}}
\includegraphics[height=2.4cm]{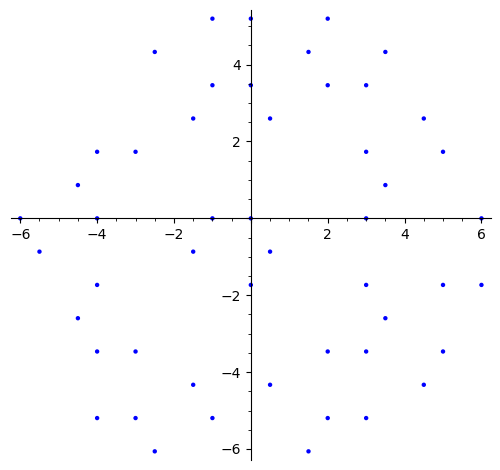} &
\includegraphics[height=2.8cm]{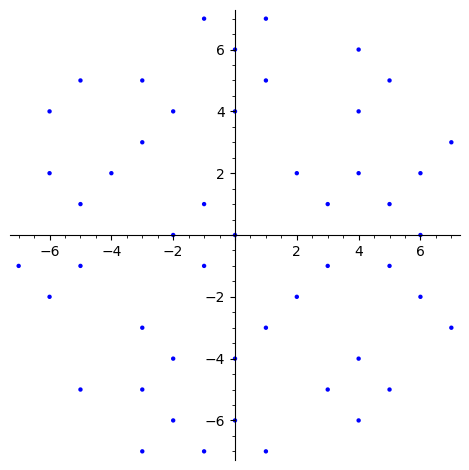} &
\includegraphics[height=4.4cm]{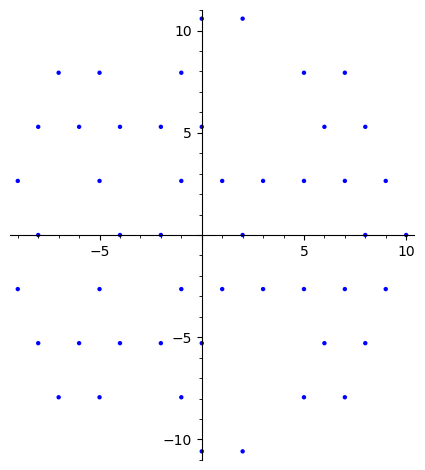}
\end{tabular}
\caption{Winning sets $\Delta=-3,-4,-7$}
\label{fig:best347}
\end{figure}
\end{center}

\begin{center}
\begin{figure}[htbp]
\centering
\begin{tabular}{m{0.33\linewidth} m{0.33\linewidth} m{0.33\linewidth}}
\includegraphics[height=4cm]{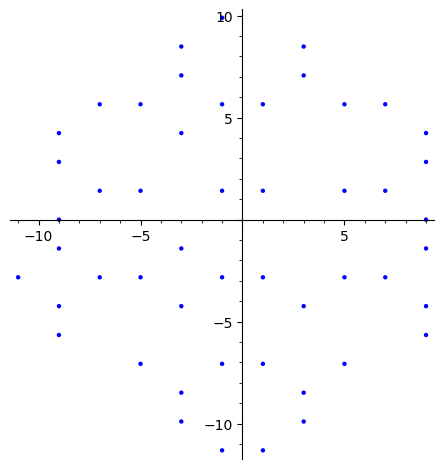} &
\includegraphics[height=4cm]{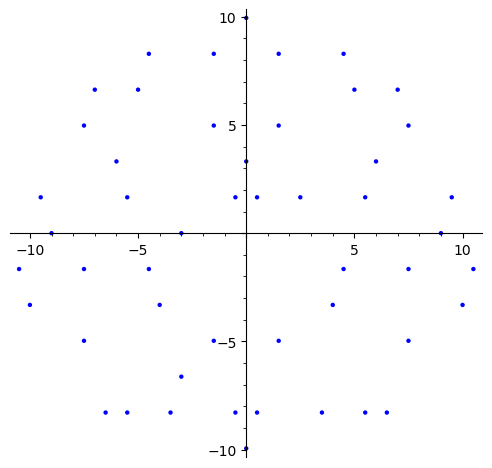} &
\includegraphics[height=4cm]{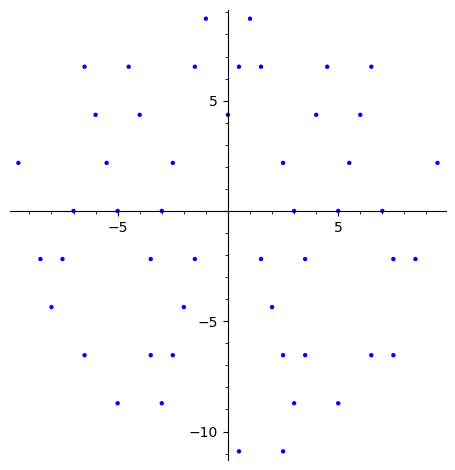}\\ 
\end{tabular}
\caption{Winning sets $\Delta=-8,-11,-19$}
\label{fig:best81119}
\end{figure}
\end{center}

\begin{center}
\begin{figure}[htbp]
\centering
\begin{tabular}{m{0.33\linewidth} m{0.33\linewidth} m{0.33\linewidth}}
\includegraphics[height=3.4cm]{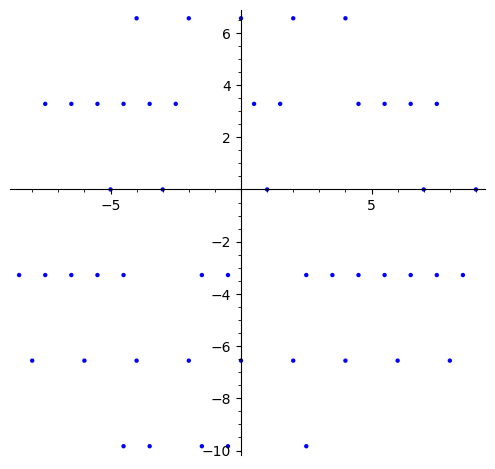} &
\includegraphics[height=3.4cm]{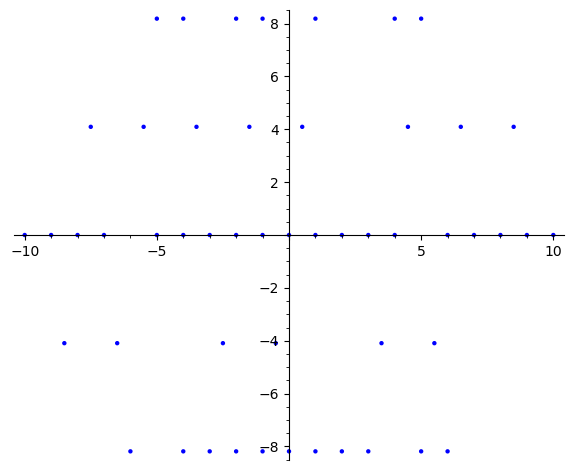} &
\includegraphics[height=4cm]{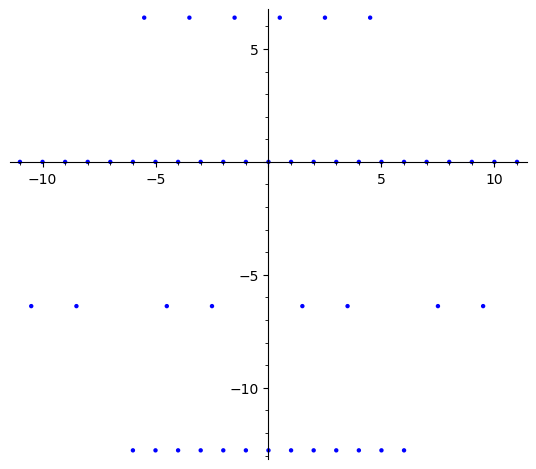}\\ 
\end{tabular}
\caption{Winning sets $\Delta=-43,-67,-163$}
\label{fig:best4367163}
\end{figure}
\end{center}

The shortest diameters found are highlighted in bold. For these sets we plot the Minkowski embeddings of the corresponding admissible sets in Figures \ref{fig:best347}, \ref{fig:best81119}, and \ref{fig:best4367163}. The specific elements of each set and the missing residue classes for each prime ideal $\mathfrak p:\mathcal N\mathfrak p<50$ are given in the appendix.

Note that the computations in bold in table \ref{tab:disc} and the associated images provide a proof of Theorem \ref{th:disc}.

The diversity in regions of bounded diameter makes an exhaustive search of the full system of residues like that of section \ref{sec:lines} more challenging. There is considerable scope for improvement of the results of this section and hopefully other researchers will be able to contribute improved sets.

\section{Extending the second Hardy-Littlewood conjecture\label{sec:ehl}}
We now turn to consideration of the extended second Hardy-Littlewood conjecture (\ref{eq:ghl2}) in our fields. We note that the ``lucky polynomial'' phenomenon observed by Euler means that the fields with larger discriminant have sparse sets of very small primes so that $\pi(0,A)$ is very small for small $A$. We anticipate that counter-examples will be easier to find in these fields and begin our search in reverse discriminant order.

For the field $\mathbb Q(\sqrt{-163})$ we note that the smallest $A_0$ such that $\pi(0,A_0)>0$ is $A_0=4$ when the disc $\mathcal N(\eta)\le 4$ captures the two associates of the inert prime 2 so that $\pi(0,4)=2$. The disc itself contains five integers -2, -1, 0, 1, and 2 and our constructions in section \ref{sec:lines} readily provide a translate of this set all of whose elements admissible (e.g. by selecting the residue class $\alpha_{163}=(1+\sqrt{-163})/2$ of the inert prime 2 for the centre element). It is then straightforward to locate a 5-tuple of primes corresponding to this admissible set $\{-2+\alpha_{163},-1+\alpha_{163},\alpha_{163},1+\alpha_{163},2+\alpha_{163}\}$.

The same admissible set construction will apply in any field where there are only two roots of unity and the primes 2, 3 and 5 are inert. We therefore handle the fields $\mathbb Q(\sqrt{-67}$ and $\mathbb Q(\sqrt{-43})$ in a similar fashion, noting that $\{-2+\alpha_{67},\ldots, 2+\alpha_{67}\}$ and $\{-2+\alpha_{43},\ldots,2+\alpha_{43}\}$ are primes in their respective fields for $\alpha_{67}=(1+\sqrt{-67})/2$ and $\alpha_{43}=(1+\sqrt{-43})/2$.

The field $\mathbb Q(\sqrt{-19})$ is only slightly different in that we must also sieve by the primes $\alpha_{19}=(1+\sqrt{-19})/2$ and $1-\alpha_{19}$ which divide 5. We must eliminate at least one prime from our admissible set due to these primes, but if we specify that the middle element is in the residue class 0 for both we still have an admissible set of four elements -2, -1, 1, 2. Again, the set of four $\{-2+\alpha_{19},-1+\alpha_{19},1+\alpha_{19},2+\alpha_{19}\}$ provides an explicit counter-example (we could also argue that we should be allowed to include the prime element $\alpha_{19}$ for a larger counter-example).

Cheap tricks such as this no longer avail us in $\mathbb Q(\sqrt{-11})$ where the prime 3 splits and is not ramified: $3=\alpha_{11}(1-\alpha_{11})$ (where $\alpha_{11}=1+\sqrt{-11})/2$). We begin by considering the value $A_0=3$ for which $\pi(0,A_0)=4$ and where the disc contains seven integers $0,\pm 1,\pm\alpha_{11},\pm(1-\alpha_{11})$, exhaustively or otherwise we soon convince ourselves that after sieving by 2, $\alpha_{11}$ and $1-\alpha_{11}$ will eliminate at least four candidates. Increasing to $A_0=4$ increases the target size of the set by 2 due to the inert prime 2 and its associate, while only increasing the initial element count by 2. We conclude that a suitable admissible set does not exist for $A_0=4$. Likewise, increasing to $A_0=5$ raises the target by 4 but only increases the number of candidates by the same amount. Resorting to the methods of section \ref{sec:discs}, attempting to greedily find an admissible set with more survivors than primes in the disc under consideration does not succeed for sets of up to 4000 candidates. This does not preclude the existence of suitable admissible sets within these candidates; it merely demonstrates the limitations of the methods in section \ref{sec:discs}.

The case $\mathbb Q(\sqrt{-8})$ is saved by the relatively large gap in prime norms between 3 and 11. This suggests targetting the value $A_0=9$ (there are no elements of norm 10) giving us 23 candidates of which at least seven must survive sieving by $\sqrt{-2}$ and $1\pm\sqrt{- 2}$. It is not hard to recover admissible sets such as $\{-1,1-\sqrt{- 2},1+\sqrt{- 2},1-2\sqrt{-2},  -3, 3, 1+2\sqrt{- 2}\}$. Explicit counter-examples such as $\{1+3\sqrt{-2}, 3+2\sqrt{-2},3+4\sqrt{-2},3+\sqrt{-2},-1+3\sqrt{-2},5+3\sqrt{-2},3+5\sqrt{-2}\}$ are then easy to find.

A similar gap in prime norms between 2 and 7 in $\mathbb Q(\sqrt{-7})$ does not quite admit an admissible set of size five in the disc $\mathbb N(\eta)\le 4$. However, expanding to the disc $\mathcal N(\eta)\le 8$ increases the number of candidates by ten while only increasing the required number of survivors by two. One quickly finds the admissible set $\{-2-2\alpha_7,-2,-2\alpha_7,0,2\alpha_7,2,2+2\alpha_7\}$. Explicit counter-examples such as $\{3+2\alpha_7,3+4\alpha_7,5+2\alpha_7,5+4\alpha_7,5+6\alpha_7,7+4\alpha_7,7+6\alpha_7\}$ are then easy to find.

In the Gaussians, the presence of additional roots of unity, combined with ramification of the smallest prime leads to a speedy early growth of candidates and a slower growth of necessary survivors. Although the minimal $A_0=2$ disc does not quite permit an admissible set of five of its nine candidates, the next value $A_0=4$ introduces four new candidates and no increase in the necessary divisors. This allows us to find the over-achieving admissible set of eight elements $\{-2i,-1-i,1-i, -2,2,-1+i,1+i, 2i\}$ of which any subset of five elements could provide a counter-example. An explicit counter-example of seven primes is then easy to find, for example $\{55+78i, 54+79i,56+79i,53+80i,57+80i,56+81i,55+82i\}$

A similar effect in the Eisenstein integers provides the minimal $A_0=3$ with 13 candidates requiring seven survivors. Again, we overachieve with admissible sets such as the nine elements $\{-1-2\omega,-2-\omega,-1-\omega,1-\omega,1,-1+\omega,\omega, 2+\omega,1+2\omega\}$, of which any subset of seven elements could provide a counter-example (here $\omega$ is a non-trivial cube root of unity). An explicit example of eight primes is $\{13+12\omega,12+13\omega, 15+13\omega,15+14\omega,13+15\omega,14+15\omega,16+15\omega,15+16\omega\}$.

\section{Lies, damned lies and heuristics}
One might ask whether the counter-examples of section \ref{sec:ehl} are artifacts of the small regions considered and whether in each field there exist counter-examples to equation (\ref{eq:ghl2}) with arbitrarily large $A$. In this section we present heuristic to suggest that for large $A$ the behaviour is commensurable with the rational case.

We note various field constants for our fields: the number of real roots $r_1=0$, the number of pairs of complex roots $r_2=1$, the class number $h=1$, the unit group is finite and hence the regulator is $R=1$, and in $\mathbb Q(\sqrt{-3})$ we have $w=6$ roots of unity, in $\mathbb Q(\sqrt{-1})$ we have $w=4$ and for the other fields we have $w=2$.

As every ideal is principal, the number of elements in the disc ${\mathcal N}(\eta)<A$ is equal to
\begin{equation}w\sum_{I:{\mathcal N}(I)<A}1\sim wA\cdot\mathrm{Res}_{s=1}\zeta_K(s)=\frac{2\pi A}{\sqrt{|\Delta|}}\label{eq:candidates}\end{equation}
by the class number formula with $\Delta$ the discriminant of our field. This is the number of candidates which we will sieve to construct admissible sets. Similarly, by the prime ideal theorem we have
\begin{equation}\pi(0,A)=w\sum_{\mathfrak p:{\mathcal N}\mathfrak p<A}1\sim \frac {wA}{\log A}.\label{eq:pit}\end{equation}
This is the target lower bound for the size of an admissible set which would (subject to the first generalised Hardy-Littlewood conjecture) provide a counter-example to (\ref{eq:ghl2}).

Now suppose that we seek an admissible set of size $z$ contained in a disc $\mathcal N(\eta)<A=o(\exp(z))$. Under the Minkowski embedding, the set corresponds to set of points in our disc modulo the lattice that is the image of the ideal $\mathfrak P:=\prod_{\mathfrak p:\mathcal N\mathfrak p<z}\mathfrak p$ intersected with points coprime to $\mathfrak P$ (the set may correspond to more than one such set of points if there exist primes where more than one residue class is absent from the admissible set). The proportion of points in our quotient by the lattice which are coprime to $\mathfrak P$ is given by 
\begin{equation}\prod_{\mathfrak p:\mathcal N\mathfrak p<z}\left(1-\frac1{\mathcal N\mathfrak p}\right)\sim\frac{e^{-\gamma}}{\frac{2\pi}{w\sqrt{|\Delta|}}\log z}\label{eq:mertens}\end{equation}
by Mertens theorem for number fields \cite{rosen1999generalization} where $\gamma$ is the usual Euler-Mascheroni constant 0.5772... Taking $z=wA/\log A$ per equation (\ref{eq:pit}) together with equations (\ref{eq:candidates}) and (\ref{eq:mertens}) we see that a random disc of area $A$ centred on $\eta\pmod{\mathfrak P}$ will on average lead to an admissible set of size asymptotic to $wA\exp(-\gamma)/\log A$ which is less than our target. However, we do not expect the size of these sets to be uniformly distributed. Large prime gap methods are able to exhibit many discs of this size where there are no admissible elements and it feels reasonable to expect some discs where the number of admissible elements exceeds the average by a factor of greater than $\exp(\gamma)$. Note that duplicating the argument in the rational field with intervals in place of discs leads to the same excess factor.

Writing $\phi_K(\mathfrak P)$ for the number of reduced residue classes for $\mathfrak P$ in our field and ignoring locality, selecting $r$ distinct residue classes at random modulo $\mathfrak P$, the number of admissible residues selected, $a$, satisfies a hyper-geometric distribution, 
$$\mathbb P(a=k)=\frac{\binom {\phi_K(\mathfrak P)}k\binom{|\mathcal N\mathfrak P|-\phi_K(\mathfrak P)}{r-k}}{\binom{|\mathcal  N\mathfrak P|}r}$$
which for small $r$ relative to $\mathcal N\mathfrak P$ tends to a Poisson distribution, with parameter equal to the expectation. Thus if we believe that the number of admissible residues in our discs matches the distribution for general subsets of residues, the admissible elements per disc will be distributed $\mathrm{Poi}(\lambda=wA\exp(-\gamma)/\log A)$. A similar combinatorial construction occurs in Gallagher's proof (subject to the first Hardy-Littlewood conjecture) that the number of primes in the interval $[x,x+\lambda\log x]$ is Poisson distributed with parameter $\lambda$ \cite{gallagher1976distribution}. In our case, the proportion of discs which exceed the average by our desired factor will be (assuming that our Poisson model is valid) will be the tail probability
$$e^{-\lambda}\sum_{j>e^{\gamma\lambda}}\frac{\lambda^j}{j!}$$
which will be dominated by the leading term
$$\frac{\lambda^{\lceil e^\gamma\lambda\rceil}}{(\lceil e^\gamma\lambda\rceil)!}e^{-\lambda}\sim\frac1{\sqrt{2\pi e^\gamma\lambda}}\exp(\lambda(e^\gamma(1-\gamma)- 1))\approx\frac1{\sqrt{2\pi e^\gamma\lambda}}\exp(-0.247\lambda).$$
The number of discs under consideration is the number of residue classes mod $\mathfrak P$ (which will be the number of possible centres for our discs). We see 
$$\log \mathcal N\mathfrak P=\sum_{\mathfrak :\mathcal N\mathfrak p<z}\log\mathcal N\mathfrak p\sim z,$$
by the prime ideal theorem, so that we expect the number of discs to be $\exp((1+o(1))wA/\log A)\approx\exp(1.781\lambda)$ which provides hope that such dense admissible sets exist.

The unreliability of such a model is underscored by the results of section \ref{sec:lines} where we construct one-dimensional admissible sets where the asymptotic excess is a factor $\Omega(\sqrt{\log z})$ rather than merely a constant.

The small explicit counter-examples of section \ref{sec:ehl} can be partly explained by an examination of the subsidiary terms in the asymptotics used above. Garcia and Simpson \cite{garcia2022unconditional} note that an explicit form of Mertens' theorem is
$$\prod_{\mathcal N\mathfrak p\le x}\left(1-\frac1{\mathcal N\mathfrak p}\right)=\frac{e^{-\gamma}}{\frac{2\pi}{w\sqrt|\Delta|}\log x}\exp(-E_K(x))$$
where $E_K(x)=F_K(x)+B_K(x)$ and
$$F_K(x)=\sum_{\mathfrak p:\mathcal N\mathfrak p>x}\left(\frac1{\mathcal N\mathfrak p}+\log\left(1-\frac1{\mathcal N\mathfrak p}\right)\right)=-\sum_{\mathfrak p:\mathcal N\mathfrak p>x}\sum_{e=2}^\infty\frac1{e\mathcal N\mathfrak p^e}$$
and
$$B_K=\frac{A_K(x)}{\log x}-\int_x^\infty\frac {A_K(t)}{t(\log t)^2}dt$$
where
$$A_K(x)=\sum_{\mathfrak p: \mathcal N\mathfrak p<x}\frac{\log\mathcal N\mathfrak p}{\mathcal N\mathfrak p}-\log x.$$
$F_K(x)$ drops off at a rate $O(1/x)$ whereas $A_K(x)$ tends to a non-zero constant. We therefore expect the ``advantage'' of a field is searching for small counter-examples to (\ref{eq:ghl2}) to be approximately quantified by more negative values of $A_K(x)/\log x$ (as well a small values of $\pi(0,A)$ which we also anticipate in our fields). We tabulate the values of this advantage measure in the table in Figure \ref{tab:disc}
\begin{center}
\begin{figure}[ht]
\centering
\begin{tabular}[h]{|c||c|c|c|c|c|}
\hline
Disc. & $A_K(4)/\log 4$ & $A_K(16)/\log 16$ & $A_K(64)/\log 64$ & $A_K(256)/\log 256$ & $\lim A_K$\\
\hline
1 ($\mathbb Q$) & -0.486 & -0.376 & -0.283 & -0.228 & -1.33 \\ 
-3 & -0.736 & -0.400 & -0.319 & -0.244 & -1.43 \\ 
-4 & -0.5 & -0.287  & -0.211 & -0.195 & -1.11 \\
-7 & -0.500 & -0.404 & -0.326 & -0.250 & -1.46 \\
-8 & -0.222  & -0.454 & -0.312 & -0.236 & -1.38\\
-11 & -0.472 & -0.300 & -0.240 & -0.212 & -1.21 \\
-19 & -1 & -0.197 & -0.168 & -0.177 & -1.00 \\
-43 & -1 & -0.487 & -0.236 & -0.181 & -1.125\\
-67 & -1 & -0.787 & -0.412 & -0.263 & -1.478 \\
-163 & -1 & -0.787 & -0.614 & -0.517 & -2.93\\
\hline
\end{tabular}
\caption{Advantage measure across fields and bounds}
\label{tab:advantage}
\end{figure}
\end{center}

\section{Conclusion}
We have demonstrated methods to construct dense admissible sets in imaginary quadratic fields of class number 1, having seen that methods used in constructing rational admissible sets do not transfer well. We have constructed especially dense admissible sets where all elements differ by rational integers. This leads to new records for the smallest rational integer difference between large primes in our fields. We have also constructed admissible sets contained within a disc of small radius, leading to new records on the smallest norm of the difference between large primes in our fields. We have used very small admissible sets to construct explicit counter-examples to the natural generalisation of the second Hardy-Littlewood conjecture in all of our field except $\mathbb Q(\sqrt{-11})$ and provided some heuristics on the expected distribution of the size of such admissible sets.

\bibliographystyle{plain}
\bibliography{admissible_gauss}

@article{renze2001gaussian,
  title={The {Gaussian} zoo},
  author={Renze, John and Wagon, Stan and Wick, Brian},
  journal={Experimental Mathematics},
  volume={10},
  number={2},
  pages={161--173},
  year={2001},
  publisher={Taylor \& Francis}
}

@article{goldston2009primes,
  title={Primes in tuples I},
  author={Goldston, Daniel A and Pintz, J{\'a}nos and Yildirim, Cem Y},
  journal={Annals of Mathematics},
  pages={819--862},
  year={2009},
  publisher={JSTOR}
}

@article{zhang2014bounded,
  title={Bounded gaps between primes},
  author={Zhang, Yitang},
  journal={Annals of mathematics},
  pages={1121--1174},
  year={2014},
  publisher={JSTOR}
}

@article{maynard2015small,
  title={Small gaps between primes},
  author={Maynard, James},
  journal={Annals of mathematics},
  pages={383--413},
  year={2015},
  publisher={JSTOR}
}

@article{polymath2014variants,
  title={Variants of the {Selberg} sieve, and bounded intervals containing many primes},
  author={Polymath, DHJ},
  journal={Research in the Mathematical sciences},
  volume={1},
  number={1},
  pages={12},
  year={2014},
  publisher={Springer}
}

@article{stadlmann2026bounded,
  title={Bounded gaps between primes},
  author={Stadlmann, Julia},
  journal={arXiv preprint arXiv:2608.31126},
  year={2026}
}

@article{song2026bounded,
  title={Bounded Gaps Between Primes: An Upper Bound of 236},
  author={Song, Zhao and Yue, Song},
  journal={Cryptology ePrint Archive},
  year={2026}
}

@article{castillo2015bounded,
  title={Bounded gaps between primes in number fields and function fields},
  author={Castillo, Abel and Hall, Chris and Lemke Oliver, Robert and Pollack, Paul and Thompson, Lola},
  journal={Proceedings of the American Mathematical Society},
  volume={143},
  number={7},
  pages={2841--2856},
  year={2015}
}

@article{vatwani2017bounded,
  title={Bounded gaps between {Gaussian} primes},
  author={Vatwani, Akshaa},
  journal={Journal of Number Theory},
  volume={171},
  pages={449--473},
  year={2017},
  publisher={Elsevier}
}

@article{hardy1923some,
  title={Some problems of ‘Partitio numerorum’; III: On the expression of a number as a sum of primes},
  author={Hardy, Godfrey H and Littlewood, John E},
  journal={Acta mathematica},
  volume={44},
  number={1},
  pages={1--70},
  year={1923},
  publisher={Springer}
}

@inproceedings{hensley1973incompatibility,
  title={On the incompatibility of two conjectures concerning primes},
  author={Hensley, Douglas and Richards, Ian},
  booktitle={Analytic Number Theory (Proc. Sympos. Pure Math., Vol. XXIV, St. Louis Univ., St. Louis, Mo. 1972)},
  pages={123--127},
  year={1973}
}

@incollection{erdos1981problems,
  author    = {Erd{\H{o}}s, Paul},
  title     = {Problems and results in number theory},
  booktitle = {Recent Progress in Analytic Number Theory, Vol. 1},
  editor    = {Halberstam, Heini and Hooley, Christopher},
  pages     = {1--13},
  publisher = {Academic Press},
  address   = {London-New York},
  year      = {1981},
  note      = {Symposium, Durham, 1979}
}

@article{clark2001dense,
  title={Dense admissible sequences},
  author={Clark, David and Jarvis, Norman},
  journal={Mathematics of computation},
  volume={70},
  number={236},
  pages={1713--1718},
  year={2001}
}

@article{schinzel1961,
author = {Schinzel, Andrzej},
journal = {Acta Arithmetica},
language = {eng},
number = {1},
pages = {1-8},
title = {Remarks on the paper "Sur certaines hypothèses concernant les nombres premiers"},
url = {http://eudml.org/doc/206986},
volume = {7},
year = {1961},
}

@article{hardy1915expression,
  title={On the expression of a number as the sum of two squares},
  author={Hardy, Godfrey Harold},
  journal={Quart. J. Math.},
  volume={46},
  pages={263--283},
  year={1915}
}

@article{rosen1999generalization,
  title={A generalization of {Mertens'} theorem},
  author={Rosen, Michael},
  journal={JOURNAL-RAMANUJAN MATHEMATICAL SOCIETY},
  volume={14},
  pages={1--20},
  year={1999},
  publisher={THE RAMANUJAN MATHEMATICAL SOCIETY}
}

@article{gallagher1976distribution,
  title={On the distribution of primes in short intervals},
  author={Gallagher, Patrick X},
  journal={Mathematika},
  volume={23},
  number={1},
  pages={4--9},
  year={1976},
  publisher={London Mathematical Society}
}

@article{garcia2022unconditional,
  title={Unconditional explicit {Mertens’} theorems for number fields and {Dedekind} zeta residue bounds},
  author={Garcia, Stephan Ramon and Lee, Ethan Simpson},
  journal={The Ramanujan Journal},
  volume={57},
  number={3},
  pages={1169--1191},
  year={2022},
  publisher={Springer}
}

\appendix

\section{Admissible set elements}

Here we provide explicit descriptions of admissible sets described in sections \ref{sec:lines} and \ref{sec:discs} for each of our fields $\mathbb Q(\sqrt{-\Delta})$.

For $\Delta=-3$ (the Eisenstein integers), we find that the 94 integers in the interval $[24943,25037]$ include 50 which are coprime to $3\cdot 7\cdot13\cdot19\cdot31\cdot37\cdot43$.

We also see that the set
\begin{eqnarray*}&&\{-\omega, 0, -1-2\omega, -1, -2-\omega, 3, 2+3\omega, 2-2\omega, 3\omega, -4\omega, 4+2\omega, 4+\omega,-2+2\omega,\\
&& -2-5\omega, 2+4\omega,1+4\omega, -4, 2-3\omega, -4-5\omega, -5-2\omega, -5-4\omega, 1-4\omega, 4+4\omega,\\
&& -3+2\omega,-4+\omega, -4-6\omega, 5+4\omega, -6-3\omega, -6-4\omega, 4+5\omega, -1-6\omega, 4-2\omega, 5\omega,\\ &&6+3\omega, 6+2\omega, -6-\omega, -6-6\omega, -6\omega, 3+6\omega,2+6\omega, -2-7\omega,3-4\omega, 6+5\omega,\\
&& 2-5\omega, 5+6\omega, 6, -6, -6-7\omega, 5-2\omega, -7-6\omega\}\end{eqnarray*}
is missing the following residue classes for the following primes ideals:

\begin{tabular}{|l|l|}
\hline
$\mathfrak p$ & Absent residue class(es) \\
\hline
$\langle\omega - 1\rangle$ & $\omega$ \\
$\langle 2\rangle$ & $\omega + 1$\\
$\langle-3\omega - 2\rangle$ & $-3\omega$\\
$\langle3\omega + 1\rangle$ & $2\omega$\\
$\langle\omega + 4\rangle$ & $-2\omega$\\
$\langle\omega - 3\rangle$ & $6\omega$\\
$\langle2\omega - 3\rangle$ & $-3\omega$\\
$\langle-2\omega - 5\rangle$ & $-3\omega$\\
$\langle5\rangle$ & $2\omega$\\
$\langle\omega + 6\rangle$ & $-14\omega, -12\omega, 10\omega$\\
$\langle\omega - 5\rangle$ & $10\omega, 15\omega$\\
$\langle3\omega - 4\rangle$ & $-5\omega, 2\omega$\\
$\langle-3\omega - 7\rangle$ & $-17\omega, -10\omega, 6\omega$\\
$\langle\omega + 7\rangle$ & $-20\omega, -11\omega, -3\omega, 8\omega, 12\omega, 17\omega, 21\omega$\\
$\langle\omega - 6\rangle$ & $-18\omega, -13\omega, -9\omega, 2\omega, 8\omega, 10\omega, 11\omega, 15\omega, 21\omega$\\
\hline
\end{tabular}

\smallskip

For $\Delta=-4$ (the Gaussian integers), we find that the 137 integers in the interval $[68917,69053]$ include 50 which are coprime to $2\cdot 5\cdot13\cdot17\cdot29\cdot37\cdot41$.

We also see that the set
\begin{eqnarray*}&&\{0, -1-i, -1+i, -2, 2-2i, 2+2i, 1-3i, 3-i, 3+i, -4i, 4i, -3-3i, -3+3i,\\ 
&&-2-4i, 4+2i, -2+4i, -4+2i, -5-i, 5-i, 5+i, 1+5i, -5+i, 4-4i, 4+4i,\\
&& -3-5i, 3-5i, -3+5i, -6i, 6, 6i, -6-2i, -2-6i, 6-2i, 6+2i, -6+2i, -7-i,\\
&&-5-5i, -1-7i, 1-7i, 5-5i, 5+5i, 1+7i, -1+7i, -5+5i, 4-6i, 4+6i,\\ 
&&-6+4i, -3-7i, 7-3i, 7+3i\}\end{eqnarray*}
is missing the following residue classes for the following primes ideals:

\begin{tabular}{|l|l|}
\hline
$\mathfrak p$ & Absent residue class(es) \\
\hline
$\langle i - 1\rangle$ & $i$\\
$\langle 2i  - 1\rangle$ & $-2i$\\ 
$\langle -2i  - 1\rangle$ & $2i$\\
$\langle 3\rangle$ & $-1$\\ 
$\langle 3i  + 2\rangle$ & $5i$\\ 
$\langle -2i  - 3\rangle$ & $3i$\\ 
$\langle i + 4\rangle$ & $-i$\\
$\langle i - 4\rangle$ & $-i$\\
$\langle 2i  - 5\rangle$ & $-9i, 11i$\\ 
$\langle -2i  - 5\rangle$ & $-2i, 8i , 9i $\\
$\langle i + 6\rangle$ & $-14i ,-10i ,7i ,9i $\\
$\langle i - 6\rangle$ & $-15i ,-7i ,-5i ,10i $\\
$\langle 5i  + 4\rangle$ & $-20i, -17i, -12i, -7i, 8i, 10i, 12i$\\ 
$\langle -4i  - 5\rangle$ & $-18i ,-14i ,-10i ,-8i ,-7i ,19i $\\
$\langle 7\rangle$ & $-3i  - 1, 3i  - 1, -2i , 2i , -i + 1, i + 1, 3i  + 1, -3i  + 2, 2, 3i  + 2,$ \\
&$-3i  + 3, -2i  + 3, 3, 3i  + 3$\\
\hline
\end{tabular}

\smallskip

For $\Delta=-7$, we find that the 127 integers in the interval $[43595,43722]$ include 50 which are coprime to $2\cdot 7\cdot11\cdot23\cdot23\cdot29\cdot43$.

We also see that the set
\begin{eqnarray*}&&\{2, -2, -2-2\alpha_7, -2\alpha_7, 2+2\alpha_7, 2\alpha_7, 2-2\alpha_7, 4+2\alpha_7, -4, -2-4\alpha_7, 2+4\alpha_7,\\ &&-6-2\alpha_7, -4-4\alpha_7, 4-2\alpha_7, 6+2\alpha_7, 4\alpha_7, -4+2\alpha_7, -6-4\alpha_7, -2+4\alpha_7, 6-2\alpha_7,\\
&&8+2\alpha_7,-8-4\alpha_7, -4-6\alpha_7, 4-4\alpha_7,8, 8+4\alpha_7, 2+6\alpha_7, -4+4\alpha_7, -8, -10-2\alpha_7,\\
&& -8-6\alpha_7, 2-6\alpha_7,8-2\alpha_7, 10+2\alpha_7, 8+6\alpha_7, -2+6\alpha_7,-8+2\alpha_7, -10-4\alpha_7,\\ 
&&6-4\alpha_7,10+4\alpha_7, -6+4\alpha_7, 10,-10-6\alpha_7, -4-8\alpha_7, 4-6\alpha_7, 10+6\alpha_7, 4+8\alpha_7,\\ 
&&-4+6\alpha_7, -2-8\alpha_7, 6+8\alpha_7\}\end{eqnarray*}
is missing the following residue classes for the following primes ideals:

\begin{tabular}{|l|l|}
\hline
$\mathfrak p$ & Absent residue class(es) \\
\hline
$\langle-\alpha_7\rangle$ & $1$\\
$\langle\alpha_7 + 1\rangle$ & $\alpha_7$\\
$\langle2\alpha_7 + 1\rangle$ & $-\alpha_7$\\
$\langle3\rangle$ & $0$\\
$\langle2\alpha_7 - 1\rangle$ & $\alpha_7$\\
$\langle-2\alpha_7 - 3\rangle$ & $-4\alpha_7$\\
$\langle2\alpha_7 - 3\rangle$ & $8\alpha_7$\\
$\langle-2\alpha_7 - 5\rangle$ & $9\alpha_7$\\
$\langle5\rangle$ & $-1$\\
$\langle4\alpha_7 + 3\rangle$ & $-4\alpha_7$\\
$\langle-4\alpha_7 - 1\rangle$ & $-12\alpha_7$\\
$\langle4\alpha_7 - 1\rangle$ & $-18\alpha_7, -6\alpha_7, 0, 16\alpha_7$\\
$\langle-4\alpha_7 - 5\rangle$ & $ -18\alpha_7, 0, 11\alpha_7, 16\alpha_7$\\
$\langle2\alpha_7 - 5\rangle$ & $-17\alpha_7, -4\alpha_7, -3\alpha_7, 0, 5\alpha_7, 15\alpha_7$\\
$\langle-2\alpha_7 - 7\rangle$ & $-4\alpha_7, 0, 6\alpha_7, 9\alpha_7, 20\alpha_7, 21\alpha_7$\\
\hline
\end{tabular}

\smallskip

For $\Delta=-8$, we find that the 169 integers in the interval $[426043,426212]$ include 50 that are coprime to $2\cdot3\cdot11\cdot17\cdot19\cdot41\cdot43$.

We also see that the set
\begin{eqnarray*}&&\{-1-2\sqrt{-2}, -1+\sqrt{-2}, -3-\sqrt{-2}, 1-2\sqrt{-2}, 1+\sqrt{-2}, -3-3\sqrt{-2},\\
&&-5-2\sqrt{-2}, 3-3\sqrt{-2}, -5+\sqrt{-2}, -3+3\sqrt{-2}, 5-2\sqrt{-2}, 5+\sqrt{-2},\\
&&-1-5\sqrt{-2}, -1+4\sqrt{-2}, 1-5\sqrt{-2}, 1+4\sqrt{-2}, -7-2\sqrt{-2}, -7+\sqrt{-2},\\ 
&&-5-5\sqrt{-2}, 7-2\sqrt{-2}, 7+\sqrt{-2}, -5+4\sqrt{-2}, -3-6\sqrt{-2}, -3+5\sqrt{-2},\\ 
&&5-5\sqrt{-2}, 5+4\sqrt{-2}, -9-\sqrt{-2}, 3-6\sqrt{-2}, 3+5\sqrt{-2}, -9, -7+4\sqrt{-2},\\ 
&&-9-3\sqrt{-2}, -9+2\sqrt{-2}, -3-7\sqrt{-2}, 9-\sqrt{-2}, 9, -3+6\sqrt{-2}, -9-4\sqrt{-2},\\
&&3-7\sqrt{-2}, 7+4\sqrt{-2}, 3+6\sqrt{-2}, -9+3\sqrt{-2}, 9-3\sqrt{-2}, 9+2\sqrt{-2},\\ 
&&-1-8\sqrt{-2}, -1+7\sqrt{-2}, -11-2\sqrt{-2}, 1-8\sqrt{-2}, 9-4\sqrt{-2}, 9+3\sqrt{-2}\}\end{eqnarray*}
is missing the following residue classes for the following primes ideals:

\begin{tabular}{|l|l|}
\hline
$\mathfrak p$ & Absent residue class(es) \\
\hline
$\langle\sqrt{-2}\rangle$ & $0$\\
$\langle-\sqrt{-2} - 1\rangle$ & $\sqrt{-2}$\\
$\langle\sqrt{-2} - 1\rangle$ & $\sqrt{-2}$\\
$\langle\sqrt{-2} + 3\rangle$ & $\sqrt{-2}$\\
$\langle\sqrt{-2} - 3\rangle$ & $\sqrt{-2}$\\
$\langle2\sqrt{-2} - 3\rangle$ & $-2\sqrt{-2}$\\
$\langle-2\sqrt{-2} - 3\rangle$ & $-2\sqrt{-2}$\\
$\langle3\sqrt{-2} - 1\rangle$ & $-7\sqrt{-2}$\\
$\langle-3\sqrt{-2} - 1\rangle$ & $9\sqrt{-2}$\\
$\langle5\rangle$ & $2\sqrt{-2}$\\
$\langle4\sqrt{-2} + 3\rangle$ & $-18\sqrt{-2}, -8\sqrt{-2}, -5\sqrt{-2}, -4\sqrt{-2}, 0, 4\sqrt{-2}, 6\sqrt{-2},$\\ 
& $17\sqrt{-2}, 20\sqrt{-2}$\\
$\langle-4\sqrt{-2} + 3\rangle$ & $-18\sqrt{-2}, -4\sqrt{-2}, 0, 3\sqrt{-2}, 4\sqrt{-2}, 6\sqrt{-2},17\sqrt{-2},$\\
& $ 20\sqrt{-2}$\\
$\langle3\sqrt{-2} + 5\rangle$ & $-20\sqrt{-2}, -6\sqrt{-2}, -3\sqrt{-2}, 2\sqrt{-2}, 5\sqrt{-2}, 8\sqrt{-2}, 15\sqrt{-2},$\\
& $19\sqrt{-2}, 20\sqrt{-2}$\\
$\langle-3\sqrt{-2} + 5\rangle$ & $-21\sqrt{-2}, -6\sqrt{-2}, -3\sqrt{-2}, -\sqrt{-2}, 2\sqrt{-2}, 5\sqrt{-2}, 8\sqrt{-2},$\\
& $19\sqrt{-2}, 20\sqrt{-2}$\\
$\langle7\rangle$ & $2\sqrt{-2} - 3, 3\sqrt{-2} - 1, -\sqrt{-2}, 0, 2\sqrt{-2}, 3\sqrt{-2}, 1, 3\sqrt{-2} + 1,$\\
& $2\sqrt{-2} + 3, 3\sqrt{-2} + 3$\\
\hline\end{tabular}

\smallskip

For $\Delta=-11$, we find that the 101 integers in the interval $[1525297,1525398]$ include 50 which are coprime to $3\cdot 5\cdot 11\cdot23\cdot31\cdot37\cdot47$. 

We also see that the set
\begin{eqnarray*}&&\{1+\alpha_{11}, \alpha_{11}, -3, 3+\alpha_{11}, 1+2\alpha_{11}, -3-3\alpha_{11}, -3\alpha_{11}, -5-\alpha_{11}, 4-\alpha_{11}, -5-2\alpha_{11}, \\
&&3-2\alpha_{11}, 3+3\alpha_{11}, 3\alpha_{11}, 6+\alpha_{11}, -5+\alpha_{11}, -5-4\alpha_{11}, 7+2\alpha_{11}, -5+2\alpha_{11},\\ &&-3-5\alpha_{11}, -2-5\alpha_{11},-8-\alpha_{11}, 7-\alpha_{11}, -6-5\alpha_{11}, 1-5\alpha_{11}, -9-3\alpha_{11}, 6-3\alpha_{11},\\ &&7+4\alpha_{11}, -3+4\alpha_{11}, 9, -9, -3-6\alpha_{11}, -8-5\alpha_{11}, 4+5\alpha_{11}, 1+5\alpha_{11}, 3-5\alpha_{11},\\
&& 9+3\alpha_{11},-6+3\alpha_{11}, 10+\alpha_{11}, -9+\alpha_{11}, -9-5\alpha_{11},4-5\alpha_{11}, 7+5\alpha_{11}, -2+5\alpha_{11},\\ &&9+4\alpha_{11},-5+4\alpha_{11}, -11-2\alpha_{11}, 9-2\alpha_{11}, -11-\alpha_{11}, 10-\alpha_{11},3+6\alpha_{11}\}\end{eqnarray*}
is missing the following residue classes for the following primes ideals:

\begin{tabular}{|l|l|}
\hline
$\mathfrak p$ & Absent residue class(es) \\
\hline
$\langle2\rangle$ & $0$\\
$\langle-\alpha_{11}\rangle$ & $-1$\\
$\langle\alpha_{11} + 1\rangle$ & $-\alpha_{11}$\\
$\langle-\alpha_{11} - 2\rangle$ & $0$\\
$\langle-\alpha_{11} + 1\rangle$ & $0$\\
$\langle2\alpha_{11} + 1\rangle$ & $5\alpha_{11}$\\
$\langle\alpha_{11} + 5\rangle$ & $-6\alpha_{11}$\\
$\langle\alpha_{11} - 4\rangle$ & $4\alpha_{11}$\\
$\langle3\alpha_{11} - 1\rangle$ & $-10\alpha_{11}, -7\alpha_{11}, 9\alpha_{11}, 13\alpha_{11}$\\
$\langle-3\alpha_{11} - 4\rangle$ & $-13\alpha_{11}, -10\alpha_{11}, 6\alpha_{11}, 11\alpha_{11}, 13\alpha_{11}$\\
$\langle3\alpha_{11} + 5\rangle$ & $-12\alpha_{11}, -8\alpha_{11}, -4\alpha_{11}, 9\alpha_{11}, 15\alpha_{11}$\\
$\langle-3\alpha_{11} + 2\rangle$ & $-18\alpha_{11}, -10\alpha_{11}, 4\alpha_{11}, 17\alpha_{11}$\\
$\langle2\alpha_{11} - 5\rangle$ & $-22\alpha_{11}, -19\alpha_{11}, -15\alpha_{11}, -9\alpha_{11}, -5\alpha_{11}, -2\alpha_{11},7\alpha_{11}, 8\alpha_{11}, $\\
&$9\alpha_{11}, 12\alpha_{11}, 14\alpha_{11}, 19\alpha_{11}$\\
$\langle-2\alpha_{11} - 7\rangle$ & $ -21\alpha_{11}, -17\alpha_{11}, -9\alpha_{11}, -7\alpha_{11}, -4\alpha_{11}, -\alpha_{11},\alpha_{11}, 7\alpha_{11},10\alpha_{11},$\\ 
& $12\alpha_{11}, 20\alpha_{11}, 22\alpha_{11}$\\
$\langle7\rangle$ & $3\alpha_{11} - 3, -\alpha_{11} - 2, 3\alpha_{11} - 2, -2\alpha_{11} - 1, -1, 3\alpha_{11} - 1, 0,$\\
& $ -3\alpha_{11} + 1, -\alpha_{11} + 1, 1, -3\alpha_{11} + 3, 3$\\
\hline\end{tabular}

\smallskip

For $\Delta=-19$, we find that the 81 integers in the interval $[663071,663052]$ include 50 which are coprime to $5\cdot7\cdot11\cdot17\cdot19\cdot23\cdot43\cdot47$.

We also see that the set
\begin{eqnarray*}&&\{-2-\alpha_{19}, 1-\alpha_{19}, 3, -3, -4-\alpha_{19}, 3-\alpha_{19}, -3-2\alpha_{19}, 1-2\alpha_{19}, 3+\alpha_{19}, -2+\alpha_{19}, \\
&&5, -5, 1+2\alpha_{19},-4-3\alpha_{19}, 1-3\alpha_{19}, 6+\alpha_{19}, -5+\alpha_{19}, -5-3\alpha_{19}, 2-3\alpha_{19}, 5+2\alpha_{19},\\ 
&&-3+2\alpha_{19}, 7, -7, -8-\alpha_{19},7-\alpha_{19}, 2+3\alpha_{19}, 3+3\alpha_{19}, 3\alpha_{19}, 7+2\alpha_{19}, -5+2\alpha_{19},\\ 
&&-5-4\alpha_{19}, 1-4\alpha_{19}, -8-3\alpha_{19},5-3\alpha_{19}, -9-\alpha_{19}, 8-\alpha_{19}, -9-2\alpha_{19}, 6+3\alpha_{19},\\ &&-3+3\alpha_{19}, -7-4\alpha_{19}, 3-4\alpha_{19}, 6-3\alpha_{19},-2-5\alpha_{19}, 3+4\alpha_{19}, 1+4\alpha_{19}, 8+3\alpha_{19},\\ 
&&-5+3\alpha_{19}, 10+\alpha_{19}, -9+\alpha_{19}, -5\alpha_{19}\}\end{eqnarray*}
is missing the following residue classes for the following primes ideals:

\begin{tabular}{|l|l|}
\hline
$\mathfrak p$ & Absent residue class(es) \\
\hline
$\langle2\rangle$ & $0$\\
$\langle-\alpha_{19}\rangle$ & $-1$\\
$\langle\alpha_{19} + 1\rangle$ & $-\alpha_{19}$\\
$\langle-\alpha_{19} - 2\rangle$ & $-3\alpha_{19}$\\
$\langle-\alpha_{19} + 1\rangle$ & $\alpha_{19}$\\
$\langle3\rangle$ & $\alpha_{19} - 1$\\
$\langle\alpha_{19} + 3\rangle$ & $-3\alpha_{19}$\\
$\langle\alpha_{19} - 2\rangle$ & $\alpha_{19}$\\
$\langle\alpha_{19} + 4\rangle$ & $-8\alpha_{19}$\\
$\langle\alpha_{19} - 3\rangle$ & $-6\alpha_{19}$\\
$\langle2\alpha_{19} + 1\rangle$ & $\alpha_{19}$\\
$\langle2\alpha_{19} - 1\rangle$ & $8\alpha_{19}, 11\alpha_{19}$\\
$\langle-2\alpha_{19} - 3\rangle$ & $-4\alpha_{19}$\\
$\langle3\alpha_{19} + 2\rangle$ & $-20\alpha_{19}, -16\alpha_{19}, -4\alpha_{19}, -3\alpha_{19}, -\alpha_{19}, \alpha_{19}, 6\alpha_{19}, 7\alpha_{19}, 8\alpha_{19}, 12\alpha_{19},$\\
&$ 15\alpha_{19}$\\
$\langle-3\alpha_{19} - 1\rangle$ & $-16\alpha_{19}, -14\alpha_{19}, -12\alpha_{19}, -11\alpha_{19}, -2\alpha_{19}, 0, 2\alpha_{19}, 4\alpha_{19}, 6\alpha_{19}, 8\alpha_{19},$\\ 
& $10\alpha_{19}, 13\alpha_{19}, 19\alpha_{19}, 20\alpha_{19}$\\
$\langle\alpha_{19} + 7\rangle$ & $ -21\alpha_{19}, -19\alpha_{19}, -17\alpha_{19}, -16\alpha_{19}, -14\alpha_{19}, -12\alpha_{19}, -8\alpha_{19}, -7\alpha_{19},$\\ 
& $0, 4\alpha_{19}, 5\alpha_{19}, 10\alpha_{19}, 15\alpha_{19}, 20\alpha_{19}, 21\alpha_{19}$\\
$\langle\alpha_{19} - 6\rangle$ & $-16\alpha_{19}, -14\alpha_{19}, -12\alpha_{19}, -11\alpha_{19}, -8\alpha_{19}, -6\alpha_{19}, -4\alpha_{19}, -3\alpha_{19},$\\
& $ -\alpha_{19},0,\alpha_{19}, 15\alpha_{19}, 17\alpha_{19}, 18\alpha_{19}, 22\alpha_{19}$\\
\hline
\end{tabular}

\smallskip

For $\Delta=-43$, we find that the 60 integers in the interval $[11306531,11306591]$ include 50 which are coprime to $11\cdot13\cdot17\cdot23\cdot31\cdot41\cdot47$.

We also see that the set
\begin{eqnarray*}&&\{
-1-\alpha_{43}, 1, -2-\alpha_{43}, 2-\alpha_{43}, -3, 3-\alpha_{43}, -5-\alpha_{43}, 4-\alpha_{43}, -1-2\alpha_{43}, 1+\alpha_{43},\\
&&2+\alpha_{43}, -3-2\alpha_{43}, 1-2\alpha_{43}, -5, -2+\alpha_{43}, -6-\alpha_{43}, 5-\alpha_{43}, -3+\alpha_{43}, -5-2\alpha_{43},\\
&& 3-2\alpha_{43},-7-\alpha_{43}, 5+\alpha_{43}, -4+\alpha_{43}, 6-\alpha_{43}, 7, 6+\alpha_{43}, -5+\alpha_{43}, -8-\alpha_{43},\\
&& 7-\alpha_{43}, -7-2\alpha_{43},5-2\alpha_{43}, 7+\alpha_{43}, -6+\alpha_{43}, -2-3\alpha_{43}, 1+2\alpha_{43}, -3-3\alpha_{43},\\
&& 3+2\alpha_{43}, -1+2\alpha_{43}, 1-3\alpha_{43},-9-\alpha_{43}, 8-\alpha_{43}, -5-3\alpha_{43}, 8+\alpha_{43}, -7+\alpha_{43},\\
&& 5+2\alpha_{43}, -3+2\alpha_{43}, 9, -6-3\alpha_{43}, -9-2\alpha_{43},7-2\alpha_{43}\}\end{eqnarray*}
is missing the following residue classes for the following primes ideals:

\begin{tabular}{|l|l|}
\hline
$\mathfrak p$ & Absent residue class(es) \\
\hline
$\langle2\rangle$ & $0$\\
$\langle3\rangle$ & $-1$\\
$\langle-\alpha_{43}\rangle$ & $0$\\
$\langle\alpha_{43} + 1\rangle$ & $-2\alpha_{43}$\\
$\langle\alpha_{43} + 2\rangle$ & $-\alpha_{43}$\\
$\langle\alpha_{43} - 1\rangle$ & $0$\\
$\langle\alpha_{43} + 3\rangle$ & $0$\\
$\langle\alpha_{43} - 2\rangle$ & $-3\alpha_{43}$\\
$\langle\alpha_{43} + 4\rangle$ & $0$\\
$\langle\alpha_{43} - 3\rangle$ & $-2\alpha_{43}, 2\alpha_{43}$\\
$\langle5\rangle$ & $-2$\\
$\langle\alpha_{43} + 5\rangle$ & $-9\alpha_{43}, -5\alpha_{43}, 5\alpha_{43}$\\
$\langle\alpha_{43} - 4\rangle$ & $-4\alpha_{43}, -2\alpha_{43}, -\alpha_{43}, 2\alpha_{43}$\\
$\langle\alpha_{43} + 6\rangle$ & $-17\alpha_{43}, -11\alpha_{43}, -4\alpha_{43}, -3\alpha_{43}, -\alpha_{43}, \alpha_{43}, 3\alpha_{43}, 10\alpha_{43}, 16\alpha_{43}, 17\alpha_{43}$\\
$\langle\alpha_{43} - 5\rangle$ & $-18\alpha_{43}, -13\alpha_{43}, -9\alpha_{43}, -5\alpha_{43}, \alpha_{43}, 5\alpha_{43}, 9\alpha_{43}, 11\alpha_{43}, 12\alpha_{43}, 18\alpha_{43},$\\
&$ 20\alpha_{43}$\\
$\langle2\alpha_{43} + 1\rangle$ & $-21\alpha_{43}, -20\alpha_{43}, -19\alpha_{43}, -10\alpha_{43}, -6\alpha_{43}, 2\alpha_{43}, 14\alpha_{43}, 18\alpha_{43}, 19\alpha_{43},$\\ 
& $20\alpha_{43}, 21\alpha_{43}$\\
$\langle2\alpha_{43} - 1\rangle$ & $-23\alpha_{43}, -22\alpha_{43}, -21\alpha_{43}, -18\alpha_{43}, -14\alpha_{43}, -2\alpha_{43}, \alpha_{43}, 6\alpha_{43}, 10\alpha_{43},$\\
&$ 16\alpha_{43}, 19\alpha_{43}, 20\alpha_{43}, 21\alpha_{43}, 22\alpha_{43}, 23\alpha_{43}$\\
$\langle-2\alpha_{43} - 3\rangle$ & $-23\alpha_{43}, -19\alpha_{43}, -15\alpha_{43}, -14\alpha_{43}, -11\alpha_{43}, -9\alpha_{43}, -8\alpha_{43}, -7\alpha_{43},$\\ 
& $-2\alpha_{43}, 6\alpha_{43}, 7\alpha_{43}, 8\alpha_{43}, 21\alpha_{43}, 22\alpha_{43}, 23\alpha_{43}$\\
$\langle7\rangle$ & $3\alpha_{43} - 3, -2, 3\alpha_{43} - 2, -3\alpha_{43} - 1, -1, 3\alpha_{43} - 1, -3\alpha_{43}, 2\alpha_{43}, 3\alpha_{43},$\\
&$ 3\alpha_{43} + 1, 2\alpha_{43} + 2, 3\alpha_{43} + 2, -3\alpha_{43} + 3, 3, 3\alpha_{43} + 3$\\
\hline\end{tabular}

\smallskip

For $\Delta=-67$, we find that the 56 integers in the interval $[2462507,2462563]$ include 50 which are coprime to $17\cdot19\cdot23\cdot29\cdot37\cdot47$.

We also see that the set
\begin{eqnarray*}&&\{
0, 1, -1, 2, -2, 3, -3, 4, -4, -1-\alpha_{67}, 1+\alpha_{67}, -1+\alpha_{67}, -3-\alpha_{67}, -5, 3-\alpha_{67},\\
&& -3+\alpha_{67}, 6, 5+\alpha_{67}, 5-\alpha_{67}, -5+\alpha_{67}, 7, -7, -7-\alpha_{67}, 7+\alpha_{67}, 8, -8, -1-2\alpha_{67}, \\
&&-2-2\alpha_{67}, -2\alpha_{67}, 2+2\alpha_{67}, 2\alpha_{67}, -3-2\alpha_{67}, 1-2\alpha_{67}, -1+2\alpha_{67}, -7+\alpha_{67},\\ &&-4-2\alpha_{67}, 2-2\alpha_{67}, 9, -9, -5-2\alpha_{67}, 5+2\alpha_{67},-3+2\alpha_{67}, -9-\alpha_{67}, 9+\alpha_{67},\\ 
&&4-2\alpha_{67}, 6+2\alpha_{67}, -4+2\alpha_{67}, 10, -10, -7-2\alpha_{67}, 5-2\alpha_{67}\}\end{eqnarray*}
is missing the following residue classes for the following primes ideals:

\begin{tabular}{|l|l|}
\hline
$\mathfrak p$ & Absent residue class(es) \\
\hline
$\langle2\rangle$ & $\alpha_{67}$\\
$\langle3\rangle$ & $-\alpha_{67} + 1$\\
$\langle-\alpha_{67}\rangle$ & $-6$\\
$\langle\alpha_{67} + 1\rangle$ & $-5\alpha_{67}$\\
$\langle\alpha_{67} + 2\rangle$ & $3\alpha_{67}$\\
$\langle\alpha_{67} - 1\rangle$ & $5\alpha_{67}$\\
$\langle\alpha_{67} + 3\rangle$ & $-4\alpha_{67}$\\
$\langle\alpha_{67} - 2\rangle$ & $-9\alpha_{67}$\\
$\langle5\rangle$ & $2\alpha_{67} - 2$\\
$\langle\alpha_{67} + 4\rangle$ & $-11\alpha_{67}, -10\alpha_{67}, -4\alpha_{67}, 11\alpha_{67}$\\
$\langle\alpha_{67} - 3\rangle$ & $-5\alpha_{67}, 5\alpha_{67}, 14\alpha_{67}$\\
$\langle\alpha_{67} + 5\rangle$ & $-18\alpha_{67}, -11\alpha_{67}, -10\alpha_{67}, -4\alpha_{67}, 3\alpha_{67}, 4\alpha_{67}, 11\alpha_{67}, 18\alpha_{67}$\\
$\langle\alpha_{67} - 4\rangle$ & $-14\alpha_{67}, -5\alpha_{67}, -4\alpha_{67}, 3\alpha_{67}, 4\alpha_{67}, 5\alpha_{67}, 13\alpha_{67},14\alpha_{67}, 15\alpha_{67}$\\
$\langle\alpha_{67} + 6\rangle$ & $-20\alpha_{67}, -19\alpha_{67}, -12\alpha_{67}, -11\alpha_{67}, -5\alpha_{67}, -4\alpha_{67}, -3\alpha_{67}, 3\alpha_{67}, 4\alpha_{67},$\\
& $ 11\alpha_{67}, 12\alpha_{67}, 18\alpha_{67}, 19\alpha_{67}, 20\alpha_{67}, 21\alpha_{67}$\\
$\langle\alpha_{67} - 5\rangle$ & $ -23\alpha_{67}, -22\alpha_{67}, -15\alpha_{67}, -14\alpha_{67}, -13\alpha_{67}, -6\alpha_{67}, -5\alpha_{67}, -4\alpha_{67},$\\
& $ \alpha_{67}, 4\alpha_{67}, 5\alpha_{67}, 12\alpha_{67}, 13\alpha_{67}, 14\alpha_{67}, 15\alpha_{67}, 21\alpha_{67}, 23\alpha_{67}$\\
$\langle7\rangle$ & $-3\alpha_{67} - 3, 3\alpha_{67} - 3, -3\alpha_{67} - 2, 3\alpha_{67} - 2, -3\alpha_{67} - 1, 3\alpha_{67} - 1, -3\alpha_{67},$\\
& $ 3\alpha_{67}, -3\alpha_{67} + 1, -\alpha_{67} + 1, 2\alpha_{67} + 1, 3\alpha_{67} + 1,-3\alpha_{67} + 2, -\alpha_{67} + 2,$\\ 
& $3\alpha_{67} + 2, -3\alpha_{67} + 3, \alpha_{67} + 3, 3\alpha_{67} + 3$\\
\hline\end{tabular}

\smallskip

Lastly, for $\Delta=-163$ we find that the 51 integers in the interval $[830591,830642]$ include 50 that are coprime to $41\cdot43\cdot47$.

We also see that the set
\begin{eqnarray*}&&\{
0, 1, -1, 1-\alpha_{163}, 2, -2, -3-\alpha_{163}, 3, -3, 3-\alpha_{163}, 4, -4, -5-\alpha_{163}, 5, -5, 6, -6, 7\\
&&, -7, 7-\alpha_{163},8, -8, -9-\alpha_{163}, 9, -1-2\alpha_{163}, -9, 1+\alpha_{163}, -2-2\alpha_{163}, -2\alpha_{163},\\
&& -1+\alpha_{163}, -3-2\alpha_{163}, 1-2\alpha_{163},3+\alpha_{163}, -4-2\alpha_{163}, 9-\alpha_{163}, 2-2\alpha_{163},\\
&& -3+\alpha_{163}, -5-2\alpha_{163}, 3-2\alpha_{163}, 10, -10, 5+\alpha_{163},-6-2\alpha_{163}, 4-2\alpha_{163},\\
&& -11-\alpha_{163}, -5+\alpha_{163}, -7-2\alpha_{163}, 5-2\alpha_{163}, 11, -11\}\end{eqnarray*}
is missing the following residue classes for the following primes ideals:

\begin{tabular}{|l|l|}
\hline
$\mathfrak p$ & Absent residue class(es) \\
\hline
$\langle2\rangle$ & $\alpha_{163}$\\
$\langle3\rangle$ & $-\alpha_{163} - 1$\\
$\langle5\rangle$ & $2\alpha_{163} - 2, 2\alpha_{163} - 1, 2\alpha_{163}, 2\alpha_{163} + 1, 2\alpha_{163} + 2$\\
$\langle-\alpha_{163}\rangle$ & $-20, -19, -18, -17, -16, -15, -14, -13, -12, 12, 13, 14, 15, 16, 17, 18,$\\
& $ 19, 20$\\
$\langle\alpha_{163} + 1\rangle$ & $-20\alpha_{163}, -19\alpha_{163}, -18\alpha_{163}, -17\alpha_{163}, -16\alpha_{163}, -15\alpha_{163}, -14\alpha_{163},$\\
& $ -13\alpha_{163}, -12\alpha_{163}, 12\alpha_{163}, 13\alpha_{163}, 14\alpha_{163}, 15\alpha_{163}, 16\alpha_{163}, 17\alpha_{163},$\\
& $ 18\alpha_{163}, 19\alpha_{163}, 20\alpha_{163}$\\
$\langle\alpha_{163} + 2\rangle$ & $-15\alpha_{163}, -14\alpha_{163}, -13\alpha_{163}, -12\alpha_{163}, -11\alpha_{163}, -10\alpha_{163}, -9\alpha_{163},$\\
& $ -8\alpha_{163}, -7\alpha_{163}, -6\alpha_{163}, 6\alpha_{163}, 7\alpha_{163}, 8\alpha_{163}, 9\alpha_{163}, 10\alpha_{163}, 11\alpha_{163},$\\
& $ 12\alpha_{163}, 13\alpha_{163}, 14\alpha_{163}, 15\alpha_{163}$\\
$\langle\alpha_{163} - 1\rangle$ & $-21\alpha_{163}, -20\alpha_{163}, -19\alpha_{163}, -18\alpha_{163}, -17\alpha_{163}, -16\alpha_{163}, -15\alpha_{163},$\\
& $ -14\alpha_{163}, -13\alpha_{163}, 12\alpha_{163}, 13\alpha_{163}, 14\alpha_{163}, 15\alpha_{163}, 16\alpha_{163}, 17\alpha_{163},$\\
& $ 18\alpha_{163}, 19\alpha_{163}, 20\alpha_{163}, 21\alpha_{163}$\\
$\langle\alpha_{163} + 3\rangle$ & $-23\alpha_{163}, -22\alpha_{163}, -21\alpha_{163}, -20\alpha_{163}, -11\alpha_{163}, -10\alpha_{163}, -9\alpha_{163},$\\
& $ -8\alpha_{163}, -7\alpha_{163}, -6\alpha_{163}, -5\alpha_{163}, 4\alpha_{163}, 5\alpha_{163}, 6\alpha_{163}, 7\alpha_{163}, 8\alpha_{163},$\\
& $ 9\alpha_{163}, 10\alpha_{163}, 11\alpha_{163}, 20\alpha_{163}, 21\alpha_{163}, 22\alpha_{163}, 23\alpha_{163}$\\
$\langle\alpha_{163} - 2\rangle$ & $-17\alpha_{163}, -16\alpha_{163}, -15\alpha_{163}, -14\alpha_{163}, -13\alpha_{163}, -12\alpha_{163}, -11\alpha_{163},$\\
& $ -10\alpha_{163}, -9\alpha_{163}, -8\alpha_{163}, -7\alpha_{163}, -6\alpha_{163}, 6\alpha_{163}, 7\alpha_{163}, 8\alpha_{163}, 9\alpha_{163},$\\
& $ 10\alpha_{163}, 11\alpha_{163}, 12\alpha_{163}, 13\alpha_{163}, 14\alpha_{163}, 15\alpha_{163}, 16\alpha_{163}$\\
$\langle7\rangle$ & $-3\alpha_{163} - 3, 2\alpha_{163} - 3, 3\alpha_{163} - 3, -3\alpha_{163} - 2, 2\alpha_{163} - 2, 3\alpha_{163} - 2,$\\
& $ -3\alpha_{163} - 1, -\alpha_{163} - 1, 2\alpha_{163} - 1, 3\alpha_{163} - 1, -3\alpha_{163}, \alpha_{163}, 2\alpha_{163}, 3\alpha_{163},$\\
& $ -3\alpha_{163} + 1, 2\alpha_{163} + 1, 3\alpha_{163} + 1, -3\alpha_{163} + 2, 2\alpha_{163} + 2, 3\alpha_{163} + 2,$\\
& $ -3\alpha_{163} + 3, 2\alpha_{163} + 3, 3\alpha_{163} + 3$\\
\hline\end{tabular}

\end{document}